\documentclass[11pt]{amsart}

\usepackage{amsfonts}
\usepackage{amssymb}
\usepackage{amscd}

\theoremstyle{plain}

\newtheorem{theorem}{Theorem}[section]
\newtheorem{proposition}[theorem]{Proposition}
\newtheorem{corollary}[theorem]{Corollary}
\newtheorem{lemma}[theorem]{Lemma}
\theoremstyle{definition}

\newtheorem{remark}[theorem]{Remark}

\newtheorem{conjecture}[theorem]{Conjecture}
\newtheorem{conjecture/question}[theorem]{Conjecture/Question}
\newtheorem{question}[theorem]{Question}
\newtheorem{remark/definition}[theorem]{Remark/Definition}
\newtheorem{terminology/notation}[theorem]{Terminology/Notation}

\newcommand{\marginlabel}[1]%
  {\mbox{}\marginpar{\raggedleft\hspace{0pt}\bfseries\sf#1}}

\def\QQ{{\mathbb Q}}
\def\PP{{\textbf P}}
\def\OO{{\mathcal O}}

\def\KK{\mathcal{K}}

\def\I{\mathcal{I}}

\def\cM{\mathcal{M}}

\def\Pic0{{\rm Pic}^0(X)}

\def\ff{\overline{\mathcal{F}}}
\def\mm{\overline{\mathcal{M}}}
\def\kk{\overline{\mathcal{K}}}
\def\K10{\overline{\mathcal{K}}}

\def\rmapdown#1{\Big\downarrow
   \rlap{$\vcenter{\hbox{$\scriptstyle#1$}}$ }}

\theoremstyle{remark}

\begin{document}

\title{\bf Effective divisors on $\mm_g$, curves on $K3$ surfaces, and
  the slope conjecture}

\author[G. Farkas]{Gavril Farkas}
\address{Department of Mathematics, Princeton University, Fine Hall,
Princeton, NJ 08544}
\email{{\tt gfarkas@math.princeton.edu}}
\thanks{Research of GF partially supported by the NSF Grant DMS-0140520}

\author[M. Popa]{Mihnea Popa}
\address{Department of Mathematics, Harvard University,
One Oxford Street, Cambridge, MA 02138}
\email{{\tt mpopa@math.harvard.edu}}
\thanks{Research of MP partially supported by the NSF Grant DMS-0200150}

\maketitle


\markboth{G. FARKAS and M. POPA}
{Effective divisors on $\mm_g$}

\section{Introduction}

In this paper we use the geometry of curves lying on $K3$ surfaces in order to obtain a number of
results about effective divisors on the moduli space of stable curves
$\mm_g$. We begin by showing two
statements on the slopes of such divisors: first that the Harris-Morrison Slope Conjecture fails to
hold on $\mm_{10}$ and second, that  in order
to compute the slope of $\mm_g$ for $g\leq 23$, one only has to look at the
coefficients of the classes $\lambda$ and $\delta_0$ in the standard expansion in terms of the generators
of the Picard group. The proofs are based on a general result providing inequalities
between the first few coefficients of effective divisors on $\mm_g$.
We then describe in detail
the divisor $\KK$ on  $\cM_{10}$ consisting of smooth sections of $K3$
surfaces, and its compactification
$\kk$ in the moduli space $\mm_{10}$ (which is the counterexample to
the conjecture mentioned above). As far as we know this is the first
intersection theoretic analysis of a
geometric subvariety on $\mm_g$ which is not of classical
Brill-Noether-Petri type, that is, a locus
of curves carrying an exceptional linear series $\mathfrak
g^r_d$. Along the way, various other results on the
Kodaira-Iitaka dimension of distinguished linear series  on certain moduli spaces of (pointed) stable
curves are obtained. We give the technical statements in what follows.
\footnote{Apart from the study of divisors on $\mm_g$,
the importance of using curves sitting on a $K3$ surface is
well-documented:
Lazarsfeld observed that there is no Brill-Noether type obstruction to
embed a curve in a $K3$ surface,
thus giving the first proof of the Brill-Noether-Petri Theorem  which
does not use degeneration to singular curves (cf. \cite{Laz}). More
recently, Voisin used the geometry
of $K3$ sections to prove the generic Green Conjecture
(cf. \cite{Voisin2}, \cite{Voisin3}). Very interesting geometry has been
developed by Mukai (cf. e.g. \cite{Mukai1}, \cite{Mukai2}), whose
results will be used below.}

On $\mm_g$ we denote by $\lambda$ the class of the Hodge line bundle, by
$\delta_0,\ldots, \delta_{[g/2]}$
the boundary classes corresponding to singular stable curves
and by $\delta:=\delta_0+\cdots +\delta_{[g/2]}$ the total boundary.
If $\mathbb E\subset \mbox{Pic}(\mm_g)\otimes \mathbb R$ is the
 cone of effective divisors, then following \cite{harris-morrison} we
define the \emph{slope function}
$s:\mathbb E\rightarrow \mathbb R\cup \{\infty\}$ by the formula
$$s(D):= \mbox{inf }\{\frac{a}{b}:a,b>0 \mbox{ such that }
a\lambda-b\delta-D\equiv
\sum_{i=0}^{[g/2]} c_i\delta_i,\mbox{ where }c_i\geq 0\}.$$
From the definition it follows that $s(D)=\infty$ unless $D\equiv
a\lambda-\sum_{i=0}^{[g/2]} b_i\delta_i$
with $a,b_i\geq 0$ for all $i$ (and it is well-known that $s(D)<\infty$ for any $D$
which is the closure of an effective divisor
on $\mathcal{M}_g$, see e.g. \cite{harris-morrison}). In this case one has that
$s(D)=a/{\mbox{min}}_{i=0}^{[g/2]}b_i.$
We denote by $s_g$ the slope of the moduli space
$\mm_g$, defined as $s_g:=\mbox{inf }\{s(D):D\in \mathbb E\}$.
The Slope Conjecture of Harris and Morrison predicts that
$s_g\geq 6+12/(g+1)$ (cf. \cite{harris-morrison}, Conjecture 0.1).
This is known to hold for $g\leq 12, g\neq 10$
(cf. \cite{harris-morrison} and \cite{tan}). Apart from the evidence coming from low genus, the conjecture was mainly based on  the large number of calculations of classes of geometric divisors on $\mm_g$ (see e.g. \cite{EH3}, \cite{Ha}).
The original paper \cite{harris-morrison} describes a number of
interesting consequences of the statement, e.g. a positive answer would imply that the Kodaira
dimension of $\cM_g$ is $-\infty$ if and only if $g\leq 22$.
Proposition \ref{nef} below provides a geometrically meaningful
refomulation of the conjecture saying that if there exists  an
effective divisor $D$ on $\mm_g$ with
slope $s(D)<6+12/(g+1)$, then $D$ has to contain the locus $\KK_g$ of $K3$ sections of genus $g$.

We consider the divisor $\KK$ on $\cM_{10}$
consisting of smooth curves lying on a $K3$ surface, and we denote by
$\K10$ its closure in $\mm_{10}$. For any $g\geq 20$, we look at the locus in $\mm_g$
of curves obtained by attaching a pointed curve of genus $g-10$ to a
curve in $\K10$ with a marked point. This gives a divisor in $\Delta_{10}
\subset \mm_g$, which we denote by $Z$.

Based on the study of curves lying on $K3$
surfaces, we establish in \S2 inequalities involving a number of coefficients
of any effective divisor coming from $\mathcal{M}_g$ in the expansion in
terms of generating classes.

\begin{theorem}\label{ineq1}
Let $D\equiv a\lambda-\sum_{i=0}^{[g/2]}b_i\delta_i$ be the closure in $\mm_g$ of
an effective divisor on $\mathcal{M}_g$.
\newline
\noindent (a) For $2\leq i\leq 9$ and $i=11$ we have
$$b_i \geq (6i+18)b_0 - (i+1)a.$$
The same formula holds for $i=10$ if $D$ does not contain the divisor
$Z\subset \Delta_{10}$.
\newline
\noindent
(b) $(g\geq 20)$ If $D$ contains $Z$, then either $b_{10}\geq 78b_0-11a$ as above, or
$$b_{10}\geq 78\cdot \frac{1155}{1262}\cdot b_0 - 11\cdot \frac {1170}{1262}\cdot a$$
(which means approximately $b_{10}\geq (71.3866...)\cdot b_0 - (10.1980...)\cdot a$).
\newline
\noindent
(c) We always have that $b_1\geq 12b_0  - a$.
\end{theorem}

\begin{corollary}\label{ineq2}
If $a/b_0\leq 71/10$, then $b_i \geq b_0$ for all $1\leq i\leq 9$.
The same conclusion holds for $i=10$ if $a/b_0\leq
\frac{88828}{12870} = 6.9019...$, and for $i=11$
if $a/b_0\leq 83/12$.
\end{corollary}

From this we can already see that the divisor $\kk \subset \mm_{10}$ provides
a counterexample to the Slope Conjecture.
Its class can be written as
$$\K10 \equiv a\lambda - b_0\delta_0  -\ldots - b_5\delta_5,$$
and by \cite{cukierman-ulmer}, Proposition 3.5, we have $a=7$ and $b_0 = 1$. In view
of Corollary \ref{ineq2}, this information is sufficent to show that the
slope of $\K10$ is smaller than expected.

\begin{corollary}\label{counterexample}
The slope of $\kk$ is equal to $a/b_0=7$, so strictly smaller than
the bound $78/11$ predicted by the Slope Conjecture. In particular $s_{10}=7$ (since
by \cite{tan} $s_{10}\geq 7$).
\end{corollary}

Theorem \ref{ineq1}  also allows us to formulate (at least up to genus
$23$, and conjecturally beyond that) the
following principle: \emph{the slope $s_g$ of $\mm_g$ is
computed by the quotient $a/b_0$ of the relevant divisors}. We have
more generally:

\begin{theorem}\label{slope}
For any $g\leq 23$, there exists  $\epsilon_g>0$ such that
for any effective divisor $D$ on $\mm_g$ with $s_g\leq s(D)\leq s_g + \epsilon_g$ we have
$s(D) = a/b_0$, i.e. $b_0\leq b_i$ for all $i\geq 1$.
\end{theorem}

\begin{conjecture}\label{slconj}
\emph{The  statement of the theorem holds in arbitrary genus.}
\end{conjecture}

We next carry out (in \S3 -- \S5) a detailed study of the compactified
$K3$ divisor $\kk$ considered above. In the
course of doing this we develop techniques (and give applications) which go beyond this example and
will hopefully also find other uses.
We prove the following result:

\begin{theorem}\label{class}
 The class of the divisor $\kk$ in $\rm{Pic}$$(\overline{\cM}_{10})$ is given by
$$\overline{\KK}\equiv 7\lambda-\delta_0-5\delta_1-9\delta_2-12\delta_3-14\delta_4-B_5\delta_5,$$
with $B_5\geq 6$.
\end{theorem}

Once more we see that the slope of $\kk$ is equal to $7$, so strictly
smaller than the bound $78/11$ predicted by the Slope Conjecture.
The first two coefficients in this expression were computed in
\cite{cukierman-ulmer}. Unfortunately we are unable to pin down $B_5$
with the methods of this paper. It would be
very surprising though if $B_5$ were not 15, as it will be clear from
the discussion below. However, the specific applications in which we
need $\overline{\KK}$ do not use the value of this coefficient.

To compute the class of $\kk$ we show that one can think of points in
$\KK$ in four different ways,
i.e. $\KK$ has four different realizations as a geometric divisor on $\cM_{10}$.

\begin{theorem}\label{equiv}
The divisor $\mathcal{K}$ can be described (set-theoretically) as any
of the following subvarieties of $\mathcal{M}_{10}$:
\begin{enumerate}
\item (By definition) The locus of curves sitting on a $K3$ surface.
\item The locus of curves $C$ with a non-surjective Wahl map $\psi_K:\wedge^2 H^0(K_C)\rightarrow H^0(3K_C)$.
\item The divisorial component of the locus of curves $C$ carrying  a semistable rank two vector
 bundle $E$ with $\wedge^2(E)=K_C$ and $h^0(E)\geq 7$.
\item The divisorial component of the locus of genus $10$ curves sitting on a quadric in an embedding $C\subset \PP^4$ with $\rm{deg}$$(C)=12$.
\end{enumerate}
\end{theorem}

We deduce this in fact by showing the equivalence of the four
descriptions over the locus of Brill-Noether general curves, whose
complement has codimension $2$ in $\cM_{10}$. Note that the
equivalence of descriptions (1) and (2) has been proved in
\cite{cukierman-ulmer}. We obtain the expression for the class of
$\kk$  as a consequence of a more general study of the
degenerations of multiplication maps for sections of line bundles
on curves. This is intimately related to characterizations (3) and
(4) above. It is important to emphasize the role of condition (3):
it shows that the divisor $\kk$ is a \emph{higher rank}
Brill-Noether divisor, more precisely one attached to rank $2$
vector bundles with canonical determinant. This was in fact the
initial motivation for our study, and as a general method it is
likely to lead to further developments. By extrapolating
descriptions (3) and (4) to other genera $g\geq 13$ one can
construct geometric divisors on $\mm_g$ containing the locus
$\KK_g$ and which we expect to provide other counterexamples to
the Slope Conjecture. For instance when $g=13$ the closure $D$ in
$\mm_{13}$ of the locus of curves $C$ of genus $13$ sitting on a
quadric in an embedding $C\subset \PP^5$,is a divisor containing
$\KK_{13}$ and we expect that $s(D)<6+12/{14}$. We plan to return
to these problems in the future.

There is a quite striking similarity between the class of
$\overline{\KK}$ in  $\mbox{Pic}(\mm_{10})$ and the class of the
Brill-Noether divisors in the next genus $g=11$: on $\mm_{11}$
there are two distinguished geometric divisors, the $6$-gonal
locus $\mm^{1}_{11,6}$ and the divisor $\mm^2_{11,9}$ of curves
with a $\mathfrak g^2_9$. These are distinct irreducible divisors
on $\mm_{11}$ having proportional classes (cf. \cite{EH3}, Theorem
1):
$$\alpha \mm^1_{11,6}\equiv \beta \mm^2_{11,9}\equiv
7\lambda-\delta_0-5\delta_1-
9\delta_2-12\delta_3-14\delta_4-15\delta_5,$$ for precisely
determined $\alpha, \beta \in \mathbb Z_{>0}$. The reason for this
resemblance is that  the coefficients of these divisors are (up to
a constant) the same as those of any other divisor whose pullback
to $\mm_{i,1}$, for a sufficient number of $i$'s less than $g$, is
a combination of generalized Brill-Noether divisors (cf. \S5 for
specific details). The similarity is surprising, since $\kk$
behaves geometrically very differently from all Brill-Noether
divisors. For instance, while it is known that all flag curves of
genus $g$ consisting of a rational spine and $g$ elliptic tails
are outside every Brill-Noether divisor, we prove that for every
$g$ they belong to the compactification $\kk_g$  in $\mm_g$ of the
$K3$ locus (cf. \S7). In the same context, we also look in \S6 at
linear systems on $\mm_g$ having the minimal slope $6+12/(g+1)$
predicted by the Slope Conjecture. Namely on $\mm_{11}$, where the
Slope Conjecture is known to hold, although there exist only the
two Brill-Noether divisors described above we show the following
(cf. Proposition \ref{big_iitaka} for a more precise statement):

\begin{proposition}
There exist effective divisors on $\mm_{11}$ of slope $7$ and having
Iitaka dimension equal to $19$.
\end{proposition}

This fact seems to contradict the hypothesis formulated in
\cite{harris-morrison} (and proved to be true for low $g$) that the Brill-Noether divisors
are essentially the only effective divisors on $\mm_g$ of slope $6+12/(g+1)$.

We conclude with a number of further applications.
As a consequence of Theorem \ref{class}, in \S7 we study the birational
nature of the moduli spaces
$\mm_{10,n}$ of stable genus $10$ curves with $n$ marked points:

\begin{theorem}\label{marked}
The Kodaira dimension of $\mm_{10,10}$ is nonnegative, while $\mm_{10,n}$ is
of general type for $n\geq 11$. On the other hand $\kappa(\mm_{10,n})=-\infty $ for $n\leq 9$.
\end{theorem}

Also, we remark already in \S2 (Remark \ref{univ_curve}) that the methods of the present paper give a
very quick proof of the fact that the Kodaira dimension of the universal curve
$\overline{\mathcal{M}}_{g,1}$ is $-\infty$ for $g\leq 15$, $g\neq 13, 14$.

\noindent \textbf{Acknowledgments.} We have benefitted from
discussions with Igor Dolgachev and
Joe Harris. We especially thank Sean Keel for suggesting a
potential connection between the Brill-Noether linear system on
$\mm_{11}$ and divisors on $\overline{\mathcal{F}}_{11}$, eventually
leading to Proposition \ref{big_iitaka}, and Rob Lazarsfeld for his help
in the proof of Theorem \ref{equiv}.

\section{Inequalities between coefficients of divisors}

In this section we give a geometric interpretation of the Slope Conjecture and
we establish
constraints on the coefficients of $\lambda,
\delta_0,\ldots, \delta_{11}$ for any effective divisor on $\mm_g$.

Given $g\geq 1$ we
consider   a   Lefschetz pencil of curves of genus $g$ lying on a
general $K3$ surface of degree $2g-2$ in $\PP^g$. This gives rise to a
curve $B$ in the moduli space $\overline{\mathcal{M}}_g$.
Note that any such Lefschetz pencil, considered as a family of curves
 over $\PP^1$, has at least
one section, since its base locus is nonempty. Such pencils $B$ fill up
the entire moduli space $\mm_g$ for $g\leq 9$ or $g=11$ (cf. \cite{Mukai1} , \cite{Mukai2}) and the divisor $\kk$ for $g=10$. For $g\geq 13$, the pencil $B$
fills up the locus $\KK_g$ of $K3$ sections of genus $g$ and  $\mbox{dim}(\KK_g)=19+g$.

\begin{lemma}\label{nos1}
We have the formulas $B\cdot \lambda =g+1$, $B\cdot \delta_0=
6g+18$ and $B\cdot \delta_j=0$ for $j\neq 0$.
\end{lemma}
\begin{proof}
The first two numbers are computed using classical formulas from
\cite{griffiths-harris}, pp. 508--509. The last assertion is
obvious since there are no reducible curves in a Lefschetz pencil.
\end{proof}

\begin{proposition}\label{nef}
Let $D$ be the closure in $\mm_g$ of an effective divisor on $\cM_g$. If $s(D)<6+12/(g+1)$, then $D$ contains the $K3$ locus $\KK_g$.
\end{proposition}
\begin{proof}
We consider as above the curve $B\subset
\mm_g$
corresponding to a Lefschetz
pencil of curves of genus $g$ on a general $K3$ surface $S$.
From
 Lemma \ref{nos1} we obtain that
 $B\cdot \delta  /B\cdot \lambda=6+12/(g+1)>s(D)$, which implies that $B \cdot D<0$ hence $B\subset D$. By varying $B$ and $S$ we obtain that $\KK_g\subset D$.
\end{proof}

\begin{remark}
 Proposition \ref{nef} shows  that the nefness
\footnote{Recall that, slightly abusively, a curve $B$ on a projective
variety $X$ is called \emph{nef}
if $B\cdot D\geq 0$ for every effective Cartier divisor $D$ on $X$.}
of $B$ would be a sufficient condition
for the Slope Conjecture to hold in genus $g$. Moreover, note (a bit
prematurely) that Theorems \ref{slope} and  \ref{ineq1} imply
that for $g\leq 23$ the Slope Conjecture in genus $g$ is
\emph{equivalent} to $B$ being a nef curve. The conjecture will fail for $g=10$ precisely
because $B\cdot \overline{\mathcal{K}}=-1$.
\end{remark}

For each $g\geq i+1$, starting with the pencil $B$ in
$\overline{\mathcal{M}}_i$ we can  construct a new
pencil $B_i$ in $\mm_g$ as
  follows: we fix a general pointed curve $(C,p)$
genus $g-i$. We then glue the curves in the
pencil $B$ with $C$ at $p$, along one of the sections corresponding
to the base points of the pencil.
We have that all such $B_i$ fill up $\Delta_i \subset \mm_g$
for $i\neq 10$, and the divisor $Z\subset \Delta_{10}$ when $i=10$.
\begin{lemma}\label{nos2}
We have $B_i\cdot \lambda = i+1$, $B_i\cdot \delta_0 = 6i+18$,
$B_i\cdot \delta_i = -1$ and $B_i\cdot \delta_j = 0$ for $j\neq 0, i$.
\end{lemma}
\begin{proof}
This follows from Lemma \ref{nos1} and from general principles,
as explained in \cite{chang-ran}, pp.271.
\end{proof}

\begin{proof}[Proof of Theorem \ref{ineq1} (a), (c)]
(a) Let us fix $2\leq i\leq 11, i\neq 10$. Since $D$ is the closure of a
divisor coming from $\mathcal{M}_g$, it cannot contain the whole boundary $\Delta_i$.
Thus we must have a pencil $B_i$ as above such that
$B_i\cdot D\geq 0$. The same thing holds true for $i=10$ if we know
that $Z$ is not contained in $D$.
But by Lemma \ref{nos2} this is precisely the statement of this part.

\noindent
(c) This is undoubtedly well known.
\end{proof}

The study of the coefficient $b_{10}$ is more involved, since in $\mathcal{M}_{10}$
the Lefschetz pencils of curves on $K3$ surfaces only fill up a divisor.
We need some preliminaries on divisors on $\mm_{g,n}$. For $0\leq i\leq g$ and $S\subset \{1,\ldots,n\}$, the boundary divisor
 $\Delta_{i:S}$  corresponds to the closure of the locus of nodal curves
 $C_1\cup C_2$, with $C_1$ smooth of genus $i$, $C_2$ smooth
 of genus $g-i$, and such that the marked points sitting on
 $C_1$ are precisely those labelled by $S$.
We also consider the divisor $\Delta_{irr}$ consisting of irreducible pointed curves
with one node. We denote by $\delta_{i:S}\in \mbox{Pic}(\overline{\mathcal{M}}_{g,n})$ the class of $\Delta_{i:S}$ and by $\delta_{irr}$ that of $\Delta_{irr}$.
 It is well known that the Hodge class $\lambda$, the boundaries $\delta_{irr}$ and $\delta_{i:S}$, and the tautological classes $\psi_i$ for $1\leq i\leq n$, freely generate
 $\mbox{Pic}(\overline{\mathcal{M}}_{g,n})$. To simplify notation, on
 $\mm_{g,1}$
we set $\delta_i:=\delta_{i:\{1\}}$ for $1\leq i \leq g-1$.
We have the following result whose proof we omit (cf. \cite{arbarello-cornalba}):

\begin{proposition}\label{attach}
If $j:\mm_{i,1}\rightarrow \mm_g$ is the map obtained by attaching a
general marked curve of genus $g-i$ to
the marked point of each genus $i$ curve, then
$$j^*(\lambda)=\lambda, \mbox{ } j^*(\delta_0)=\delta_0, \mbox{
}j^*(\delta_i)=-\psi+\delta_{2i-g}
\mbox{ and } j^*(\delta_k)=\delta_{i-k}+\delta_{i+k-g} \mbox{ for }  k\neq 0,i.$$
(Here we make the convention $\delta_k:=0$, for $k<0$.)
\end{proposition}

 Let  $\pi:\overline{\mathcal{M}}_{g,1}\rightarrow \mm_g$
be the forgetful morphism. We will need the following results which follow essentially from \cite{arbarello-cornalba}, \S1:
\begin{lemma}\label{push}
One has the following relations:
$$\pi_*(\lambda^2)=\pi_*(\lambda \cdot \delta_i)=\pi_*(\delta_{0}\cdot \delta_i)=0
{\rm~ for ~all~ } i=0,\ldots, g-1,\mbox{ }
\pi_*(\psi^2)=12\lambda-\delta,$$
$$\pi_*(\lambda\cdot \psi)=(2g-2)\lambda,\mbox{ } \pi_*(\psi\cdot \delta_0)=(2g-2)\delta_0,
\mbox{ } \pi_*(\psi\cdot \delta_i)=(2i-1)\delta_i {\rm~ for~} i\geq 1,$$
$$\pi_*(\delta_i^2)=-\delta_i {\rm~for~} 1\leq i\leq g-1,\mbox{ }\pi_*(\delta_i\cdot \delta_{g-i})=\delta_i,
{\rm~for~ } 1\leq i<g/2, {\rm ~and~ } $$
$$\pi_*(\delta_i\cdot \delta_j)=0
{\rm~ for ~ all ~} i,j\geq 0 {\rm~ with ~} i\neq j, g-j.$$
\end{lemma}

\noindent
We consider the \emph{Weierstrass divisor} in $\mathcal{M}_{g,1}$
$$\mathcal{W}:=\{[C,p]\in \mathcal{M}_{g,1}:p\in C \mbox{ is a
Weierstrass point}\},$$
and denote by $\overline{\mathcal{W}}$ its closure in $\overline{\mathcal{M}}_{g,1}$.
Its class has been computed by Cukierman \cite{cukierman}:
$$\overline{\mathcal{W}}\equiv -\lambda+\frac{g(g+1)}{2}\psi-\sum_{i=1}^{g-1}
{g-i+1 \choose 2}\delta_i.$$

\begin{proposition}
If $\pi:\overline{\mathcal{M}}_{g,1}\rightarrow \overline{\mathcal{M}}_g$
is the forgetful morphism, then $\pi_*(\overline{\mathcal{W}}^2)$ is an effective
divisor class on $\overline{\mathcal{M}}_g$.
\end{proposition}

\begin{proof} From the previous Lemma we have that
$$\pi_*(\overline{\mathcal{W}}^2)
\equiv a\lambda-\sum_{i=0}^{[g/2]}b_i\delta_i,$$
where $a=g(g+1)(3g^2+g+2), b_0=g^2(g+1)^2/4$ while for $1\leq i<g/2$
we have $b_i=i(g-i)(g^3+3g^2+g-1)$.
When $g$ is even $b_{g/2}=(8g^5+28g^3+33g^4+4g^2)/64.$
On the other hand we have expressions for the classes of distinguished
geometric divisors on $\overline{\mathcal{M}}_g$: when $g+1$
is composite, by looking at Brill-Noether divisors one sees
that the class
$$(g+3)\lambda-\frac{g+1}{6}\delta_0-\sum_{i=1}^{[g/2]}i(g-i)\delta_i$$
is effective (cf. \cite{EH3}, Theorem 1). When $g+1$ is
prime one has to use instead the class of
the Petri divisor, which we do not reproduce here (cf. \cite{EH3}, Theorem 2).
In either case, by comparing the coefficients $a,b_i$ above with those of
these explicit effective classes, one obtains an
effective representative for $\pi_*(\overline{\mathcal{W}}^2)$.
For instance when $g+1$ is composite it is enough to check that $b_0/a\leq (g+1)/(6g+18)$
and that $b_i/a\leq i(g-i)/(g+3)$ for $i=1,\ldots,[g/2]$, which is immediate.
\end{proof}

\begin{corollary}\label{wei}
Let $D$ be any effective divisor class on $\overline{\mathcal{M}}_{g,1}$.
Then $\pi_*(\overline{\mathcal{W}}\cdot D)$ is an effective class on
$\overline{\mathcal{M}}_g$.
\end{corollary}

\begin{proof}[Proof of Theorem \ref{ineq1} (b)]
Assume that $b_{10}<78b_0-11a$.
We consider the map
$$j:\overline{\mathcal{M}}_{10,1}\longrightarrow \mm_g$$
obtained by attaching a fixed general pointed
curve of genus $g-10$ to any curve of genus $10$ with a marked point.
Our assumption says that $B \cdot j^*(D)<0$,
where $B\subset \overline{\mathcal{M}}_{10,1}$ denotes the
curve in the moduli space coming from a Lefschetz pencil of
pointed curves of genus $10$ on a general $K3$ surface.
We can write $j^*(D)=m\pi^*(\overline{\mathcal{K}})+E,$ where $m \in \mathbb Z$ and  $E$ is an effective
divisor not containing the irreducible divisor $\pi^*(\overline{\mathcal{K}})$.
From Proposition \ref{attach} we have that
$j^*(\lambda)=\lambda, j^*(\delta_0)=\delta_0$ and $j^*(\delta_{10})=-\psi$,
while $\overline{\mathcal{K}}\equiv 7\lambda-\delta_0-\cdots$, hence
$$E\equiv (a-7m)\lambda+b_{10}\psi+(m-b_0)\delta_0+
{\rm~(other~ boundaries)}.$$ Moreover from Lemma \ref{nos1} we obtain that
$R\cdot \pi^*(\kk)=\pi_*(R)\cdot \kk=-1$, which yields the inequality
\begin{equation}\label{mm}
m\geq -B\cdot j^*(D) =-11a+78b_0-b_{10}>0.
\end{equation} Since $E$ is an effective class on $\mm_{10.1}$, from Corollary \ref{wei}
it follows that $\pi_*(\overline{\mathcal{W}}\cdot E)$ is an
effective class on $\overline{\mathcal{M}}_{10}$. An easy calculation
using
Lemma \ref{push} shows that
$$\pi_*(\overline{\mathcal{W}}\cdot E)\equiv \bigl(642b_{10}+990(a-7m)\bigr)
\lambda-55\bigl(b_{10} +18(b_0-m)\bigr)\delta_0-\cdots.$$ We now
use that for every effective divisor on $\mm_g$ the coefficient
$a$ of $\lambda$ is nonnegative \footnote{For the reader's
convenience we recall that this follows from the fact that $B\cdot
\lambda>0$ for any irreducible curve $B\subset \mm_g$ such that
$B\cap \mathcal{M}_g\neq \emptyset$, while there is always a
complete curve in $\mathcal{M}_g$ passing through a general point
(outside of the given divisor). The inequality $B\cdot \lambda>0$
is a consequence of writing $\lambda=(\lambda-\epsilon
\delta)+\epsilon \delta$ where $0<\epsilon <1/{11}$ and of using
that $\lambda-\epsilon \delta$ is ample on $\mm_g$ (cf.
\cite{CH}).}. From the previous formula we get an inequality which
combined with (\ref{mm}) yields, after a simple computation
$$b_{10}\geq 78\cdot \frac{1155}{1262}\cdot b_0 - 11\cdot \frac
{1170}{1262}\cdot a
~(~=(71.3866...)\cdot b_0 - (10.1980...)\cdot a).$$
\end{proof}

\begin{question}
Do we always have the inequality $b_{10}\geq 78b_0-11a$?
\end{question}

\subsection*{\bf Slopes of divisors and further remarks.}
The inequalities established in Theorem \ref{ineq1} allow
us to show that, at least up to genus $23$, if the slope of an
effective divisor is sufficiently small, then it is computed by the
ratio $a/b_0$.

\begin{proof}[Proof of Theorem \ref{slope}]
When $g$ is such that $g+1$ is composite, we have that $s_g\leq
6+12/(g+1)$
(this being the slope of any Brill-Noether divisor). When $g$ is even, one has the estimate
$s_g\leq \frac{2(3g^2+13g+2)}{g(g+2)}$ (this being the slope of the
Petri divisor,
cf. \cite{EH3}, Theorem 2).
It follows that  for any $g\leq 23$ there exists
a positive number $\epsilon_g$ such that
$$s_g+\epsilon_g\leq 6+ \frac{11}{i+1} {\rm~for~all~} i\leq [g/2](\leq 11).$$

Assume first that $2\leq i\leq 9$ or $i=11$. Then by Theorem
\ref{ineq1}(a)
we know that $b_i \geq (6i+18)b_0 - (i+1)a$, and so certainly $b_i\geq b_0$
if $s(D)\leq 6+ \frac{11}{i+1}$.
For $i=10$ we apply \ref{ineq1}(b): if the inequality $b_{10}\geq
78b_0 -11a$ holds, then the argument is identical. If not, we have the
inequality $b_{10}\geq (71.3866...)\cdot b_0 - (10.1980...)\cdot a$.
Thus $b_{10}\geq b_0$
as soon as the inequality $a/b_0\leq 6.9$ is satisfied. But for
$g\geq 20$ the inequality $s_g < 6.9$ holds, based on the same estimates as above.
For $i=1$, the condition
is even weaker because of the formula $b_1\geq 12b_0-a$ in \ref{ineq1}(c).
Thus the slope of $D$ is computed by $a/b_0$.
\end{proof}

\begin{remark}\label{univ_curve}
An amusing consequence of Corollary \ref{wei} is that
the Kodaira dimension of the universal curve $\mathcal{M}_{g,1}$
is $-\infty$ for all $g\leq 15$, with $g\neq 13,14$ (this
can be proved directly when $g\leq 11$).
Indeed, if we assume that some multiple of the canonical class $K_{\overline{\mathcal{M}}_{g,1}}\equiv
13\lambda+\psi-3(\delta_1+\delta_{g-1})-2\sum_{i=2}^{g-2}\delta_i$  is effective on
$\overline{\mathcal{M}}_{g,1}$, then by Corollary \ref{wei}
the same multiple of the class $D:=\pi_*(K_{\overline{\mathcal{M}}_{g,1}}\cdot \overline{\mathcal{W}})$
is effective on $\overline{\mathcal{M}}_g$.
It turns out that
$s(D)=\frac{2(13g^3+6g^2-9g+2)}{g(g+1)(4g+3)}$, and from the
definition of the slope of $\mathcal{M}_g$ we have
that $s(D)\geq s_g$. But this contradicts  the estimates on $s_g$
from \cite{tan} and \cite{chang-ran}.
\end{remark}

\section{The four incarnations of the divisor $\KK$}

This section is devoted to the proof of Theorem \ref{equiv}.
We start by reviewing some notation.
Let $\mathcal{F}_g$ be the moduli space of canonically polarized
$K3$ surfaces $(S,H)$ of genus $g$. We consider the $\PP^g$-bundle
$\mathcal{P}_g=\{(S,C):C\in |H|\}$ over $\mathcal{F}_g$ which
comes equipped with a natural rational map
$\phi_g:\mathcal{P}_g - ->\mathcal{M}_g$. By Mukai's results \cite{Mukai1}
 this
map is dominant if and only if $2\leq g\leq 9$ or $g=11$. Moreover, $\phi_g$ is generically finite if and only if $g=11$ or $g\geq 13$. For $g=10$ the map $\phi_g$ has fibre dimension $3$, whereas the fibre dimension of $\phi_{12}$ is $1$. The non-finiteness of $\phi_g$ is due to the existence of Fano threefolds $X_{2g-2}\subset \PP^{g+1}$ of genus $g=10,12$.

We denote by $\KK_g$ the closure in $\cM_g$ of the image of $\phi_g$. Thus $\KK=\KK_{10}$ is an irreducible divisor on $\cM_{10}$. By $\overline{\KK}_g$ we shall denote the Deligne-Mumford closure of $\KK_g$ in $\overline{\cM}_{g}$.

For a smooth curve $C$ and a line bundle $L$ on $C$ we consider the \emph{Wahl map} $\psi_L:\wedge^2 H^0(L)\rightarrow H^0(K_C\otimes L^{\otimes ^2})$ sending a pencil in the linear system $|L|$ to its ramification divisor, that is,
$$\psi_L(\sigma\wedge \tau):=\sigma\ d\tau-\tau\ d\sigma, \mbox{ for } \sigma,\tau\in H^0(L).$$

J. Wahl proved that if a smooth curve $C$ sits on a $K3$ surface then $\psi_{K}$ is not surjective (cf. \cite{Wahl}). This is the only known intrinsic characterization of curves sitting on a $K3$ surface. Later, Cukierman and Ulmer showed that for $g=10$ the converse of Wahl's theorem also holds (cf. \cite{cukierman-ulmer}). It is believed that for Brill-Noether-Petri general curves $C$ of genus $g\geq 13$ the non-surjectivity of $\psi_K$ is equivalent to the existence of a $K3$ surface containing $C$.  From now on we fix $g=10$ and let $C$ be a Brill-Noether general curve of genus $10$. Then $C$ carries finitely many base point free pencils $\mathfrak g^1_6$. The dual linear series $\mathfrak g^4_{12}=|K_C-\mathfrak g^1_6|$ yield embeddings $C\subset \PP^4$ with $\mbox{deg}(C)=12$. We  show that we can interpret points in the divisor $\KK$ in four geometrically meaningful ways:

\noindent
\emph{ Proof of Theorem \ref{equiv}}. \ Since being Brill-Noether special is a condition of codimension $2$ on $\mm_{10}$, the general point of each component of the divisors on $\mathcal{M}_{10}$ defined by $(1)-(4)$ will correspond to a Brill-Noether general curve. The fact that $(1)$ and $(2)$ are equivalent is  the main result of \cite{cukierman-ulmer}. The implications $(1)\Rightarrow (3)$ and $(2)\Rightarrow (4)$ were proved by Voisin (cf. \cite{Voisin1}, Proposition 3.2). We are left with showing that $(3)\Rightarrow (4)$ and that $(4)\Rightarrow (2)$.

\noindent {\textbf{(3)$\Rightarrow$(4).}} Let $C$ be a Brill-Noether general curve of genus $10$ and $E\in SU(2,K_C)$ a vector bundle with $h^0(E)\geq 7$. We denote by $\OO_C(D)$ the maximal sub-line bundle of $E$. We have the exact sequence:
\begin{equation}\label{ext}
 0\longrightarrow \OO_C(D)\longrightarrow E \longrightarrow K_C(-D)\longrightarrow 0.
\end{equation}
Using \cite{pr} Lemma 3.9 we can assume that $h^0(D)\geq 2$, hence $\mbox{deg}(D)\geq 6$ ($C$ is not $5$-gonal). From the exact sequence and from Riemann-Roch we get that
$$7\leq h^0(E)\leq h^0(D)+h^0(K_C-D)=2h^0(D)-\mbox{deg}(D)+9.$$
Again, since the curve is general, one can easily see that this can happen only when
$\mbox{deg}(D)=6$ and $h^0(E)=h^0(D)+h^0(K_C-D)$, that is,  when the cohomology coboundary map  $\delta_E\in (H^0(K_C-D)^{\vee})^{\otimes ^2}$ associated to (\ref{ext}) is zero. One can find a semi-stable bundle $E\in \mbox{Ext}^1\bigl(K_C(-D),\OO_C(D)\bigr)=H^0(2K_C-2D)^{\vee}$ with $\delta_E=0$, precisely when the multiplication map
$$\mu(D):\mbox{Sym}^2H^0(K_C-D)\rightarrow H^0(2K_C-2D)$$
is not surjective. Since both vector spaces have dimension $15$, this is equivalent to saying that
the image of the embedding $C\stackrel{|K_C-D|}\longrightarrow \PP^4$ sits on a quadric.

\noindent {\textbf{(4)$\Rightarrow$(2).}} This was shown to us by R. Lazarsfeld. Again, we fix $C$ a Brill-Noether general curve of genus $10$ and $A\in W^1_6(C)$ such that the multiplication map $\mu(A):\mbox{Sym}^2 H^0(K_C-A)\rightarrow H^0(2K_C-2A)$ is not surjective. We get from Riemann-Roch that $\mbox{dim Sym}^2H^0(K_C-A)=\mbox{dim } H^0(2K_C-2A)=15$, which implies that $\mu(A)$ is not injective either. Let $s_1,s_2\in H^0(A)$ be independent sections. One has the following commutative diagram:
$$\begin{array}{ccc}

     \mbox{Sym}^2 H^0(K_C-A)&    \stackrel{\mu(A)}\longrightarrow  & H^0(2K_C-2A) \\
      \rmapdown{\alpha} & \; & \rmapdown{\beta}  \\
    \wedge^2 H^0(K_C) &     \stackrel{\psi_K}\longrightarrow & H^0(C,3K_C) \\
\end{array}$$

Here $\alpha(\epsilon \cdot \eta):=\epsilon s_1 \wedge \eta s_2-\epsilon s_2\wedge \eta s_1$  for $\epsilon, \eta\in H^0(K_C-A)$, while $\beta(u):=(us_1)\ ds_2-(us_2)\ ds_1,$ for $u\in H^0(2K_C-2A).$ To conclude that $\psi_K$ is not surjective (or equivalently injective, since $\mbox{dim}\wedge^2H^0(K_C)=\mbox{dim }H^0(3K_C)=45$), it suffices to prove that $\alpha$ is injective.

To achieve this we define the map
$$
m:=(Pe(A)\otimes Id_{H^0(K_C)})\circ (Id_{H^0(K_C-A)}\otimes f):H^0(K_C-A)^{\otimes ^2}\rightarrow H^0(K_C)^{\otimes ^2},$$
where $Pe(A):H^0(A)\otimes H^0(K_C-A)\rightarrow H^0(K_C)$ is the Petri map, while the map $f:H^0(K_C-A)\rightarrow H^0(A)\otimes H^0(K_C)$ is defined by
$$f(\eta):=s_1\otimes (\eta s_2)-s_2 \otimes (\eta s_1).$$
\noindent One can check that
$m(\epsilon \otimes \eta)=\epsilon s_1\otimes \eta s_2-\epsilon s_2 \otimes \eta s_1,$
where $\epsilon, \eta \in H^0(K_C-A)$, which shows that
$$m(\epsilon \otimes \eta +\eta \otimes \epsilon)=\epsilon s_1\wedge \eta s_2-\epsilon s_2 \wedge \eta s_1.$$
Thus we get that $m\bigl(\mbox{Sym}^2H^0(K_C-A)\bigr)\subset \wedge^2 H^0(K_C)$, so $\alpha=m_{|_{\mbox{Sym}^2 H^0(K_C-A)}}$.
It is easy to prove that $f$ is always injective, thus when $Pe(A)$ is injective, it follows that $m$ is injective as well.

We are left with the case when $Pe(A)$ is not injective (this is a divisorial condition on $\mathcal{M}_{10}$ and cannot be ruled out by a dimension count). From the base-point-free pencil trick, $\mbox{Ker}(Pe(A))$ is one-dimensional and is generated by the element $s_1\otimes (us_2)-s_2\otimes (us_1)$, where $0\neq u\in H^0(K_C-2A)$. Then one checks that
$$\mbox{Im}\bigl(Id_{H^0(K_C-A)}\otimes f\bigr)\cap \bigl(\mbox{Ker}(Pe(A))\otimes Id_{H^0(K_C)}\bigr) =0$$
$\bigl($inside $H^0(K_C-A)\otimes H^0(A)\otimes H^0(K_C)\bigr),$
which shows that $m$ is injective in this case as well. \hfill $\Box$

\begin{remark}
Theorem \ref{equiv} establishes only a set-theoretic equality between the loci in $\cM_{10}$ described by the conditions $(1)-(4)$, that is, we do not necessarily have equalities between the appropriate cycles with the  multiplicities coming from the natural scheme structures. For instance if $\mathcal{W}$ denotes the locus of smooth curves with a non-surjective Wahl map viewed as a determinantal variety, then we have the equality of divisors $\mathcal{W}=4\mathcal{K}$ (cf. \cite{cukierman-ulmer}). This can be interpreted as saying that the corank of the Wahl map $\psi_K$ is equal to $4$ for a general $[C]\in \KK$.
\end{remark}
\begin{remark} If we use description (3) of $\mathcal{K}$, it turns out that for a general $[C]\in \KK$, the rank two vector bundle $E\in SU_2(K_C)$ with $h^0(E)\geq 7$
is unique. More precisely, if $C$ is a section of a $K3$ section $S$ then $E=\mathcal{E}_{|_C}$, where $\mathcal{E}$ is a rank
two vector bundle on $S$, which is an elementary transformation of the trivial bundle $H^0(C,A)\otimes \OO_S$ along $C$. Remarkably, the vector bundle constructed in this way does not depend on the choice of the pencil $A\in W^1_6(C)$ (cf. \cite{Voisin1}).
\end{remark}

\section{Limit linear series and degeneration of multiplication maps}

The aim of this section is to understand the following situation: suppose $\{L_b\}_{b\in B^*}$ and $\{M_b\}_{b\in B^*}$ are two families of line bundles over a $1$-dimensional family of smooth curves $\{X_b\}_{b\in B^*}$, where $B^*=B-\{b_0\}$ and $b_0\in B$. We ask what happens to the multiplication map
$$\mu_b:H^0(X_b,L_b)\otimes H^0(X_b,M_b)\longrightarrow H^0(X_b, L_b\otimes M_b)$$
as $X_b$ degenerates to a singular curve of compact type $X_{b_0}$? The answer will be given in terms of limit linear series. Everything in this section is contained (at least implicitly) in \cite{EH1} and  \cite{EH4}.

First we recall a few definitions. We fix a smooth curve $C$ and a point $p\in C$. If $l=(L,V)$ is a linear series $\mathfrak g^r_d$ on $C$ with $L\in \mbox{Pic}^d(C)$ and $V\subset H^0(L)$, then by ordering the finite set $\{\mbox{ord}_p(\sigma)\}_{\sigma \in V}$ we obtain the \emph{vanishing sequence} of $l$ at $p$
$$a^l(p):0\leq a_0^l(p)<\ldots <a_r^l(p)\leq d.$$
The \emph{weight} of $p$ with respect to $l$ is defined as $w^l(p):=\sum_{i=0}^r(a_i^l(p)-i).$

For line bundles  $L$ and $M$ on $C$ and for an element $\rho \in H^0(L)\otimes H^0(M)$, we write that $\mbox{ord}_p(\sigma)\geq k$ if $\rho$ is in the linear span of the elements of the form $\sigma \otimes \tau$, with $\sigma\in H^0(L), \tau \in H^0(M)$ and such that $\mbox{ord}_p(\sigma)+\mbox{ord}_p(\tau)\geq k$.

Let $\mu_{L,M}:H^0(L)\otimes H^0(M)\longrightarrow H^0(L\otimes M)$ be the multiplication map and $\rho \in \mbox{Ker}(\mu_{L,M})$. We shall often use the following simple fact: if $\{\sigma_i\}\subset H^0(L)$ and $\{\tau_j\}\subset H^0(M)$ are bases of the spaces of global sections \emph{adapted to the point}
$p\in C$, that is, satisfying the conditions
$$\mbox{ord}_p(\sigma_i)=a_i^L(p) \mbox{ and }
\mbox{ord}_p(\tau_j)=a_j^M(p) \mbox{ for all relevant } i\mbox{ and }j,$$
then there are distinct pairs of integers $(i_1,j_1)\neq (i_2,j_2)$ such that
$$\mbox{ord}_p(\rho)=\mbox{ord}_p(\sigma_{i_1})+\mbox{ord}_p(\tau_{j_1})=
\mbox{ord}_p(\sigma_{i_2})+\mbox{ord}_p(\tau_{j_2}).$$

Suppose now that $\pi:X\rightarrow B$ is a family of genus $g$ curves over $B={\rm Spec}(R)$, with $R$ being a complete DVR with local parameter $t$, and let $0, \eta$ denote the special and the generic point of $B$ respectively. Assume furthermore that $X_{\eta}$ is smooth and that $X_0$ is singular but of compact type. If $L_{\eta}$ is a line bundle on $X_{\eta}$ then, as explained in \cite{EH1}, there is a canonical way to associate to each component $Y$ of $X_0$ a line bundle $L^Y$ on $X$ such that $\mbox{deg}_Z(L^Y_{|_Z})=0$ for every component $Z$ of $X_0$ different from $Y$.
We set $L_Y:=L^Y_{|_Y}$ which is a line bundle on the smooth curve $Y$.

We fix $\sigma \in \pi_*L_{\eta}$ a section on the generic fibre. We denote by $\alpha$ the smallest integer such that $t^{\alpha}\sigma \in \pi_*L^Y$, that is, $t^{\alpha}\sigma \in \pi_*L^Y-t\pi_*L^Y$. Then we set
$$\sigma^Y:=t^{\alpha} \sigma \in \pi_*L^Y \mbox{ and } \sigma_Y:=\sigma^Y_{|_Y}\in H^0(Y,L_Y).$$
For a different component $Z$ of the special fibre $X_0$ meeting $Y$ at a point $p$, we define similarly $L^Z,L_Z,\sigma^Z$ and $\sigma_Z$. We have the following compatibility relation between $\sigma_Y$ and $\sigma_Z$ (cf. \cite{EH1}, Proposition 2.2):
\begin{lemma}\label{limit}
If we write $\sigma^Z=t^{\beta}\sigma^Y \in \pi_*L^Z$ for a unique integer $\beta$, then
\begin{center}
$\rm{deg}$$(L_Y)-\rm{ord}$$_p(\sigma_Y)\leq \beta \leq \rm{ord}$$_p(\sigma_Z).$
\end{center}
\end{lemma}
An immediate consequence of this is the inequality $$\mbox{ord}_p(\sigma_Y)+\mbox{ord}_p(\sigma_Z)\geq \mbox{deg}(L_Y)=\mbox{deg}(L_Z).$$

Assume from now on that we have two line bundles $L_{\eta}$ and $M_{\eta}$ on $X_{\eta}$ and we choose an element $\rho \in H^0(X_\eta,L_\eta)\otimes _{R_\eta } H^0(X_\eta, M_\eta).$ If $Y$ and $Z$ are components of $X_0$ meeting at $p$ as above, we define $\rho^Y:=t^{\gamma}\rho \in H^0(X, L^Y)\otimes _R H^0(X,M^Y)$, where $\gamma$ is the minimal integer with this property. We have a similar definition for $\rho^Z\in H^0(X,L^Z)\otimes _R H^0(X,M^Z)$.  Between the sections $\rho^Y$ and $\rho^Z$ there is a relation
$$\rho^Z=t^{\alpha} \rho^Y$$
 for a uniquely determined integer $\alpha$. To determine $\alpha$ we proceed as follows: we choose bases of sections $\{\sigma_i=\sigma^Y_i\}$ for $H^0(X,L^Y)$ and $\{\tau_j=\tau^Y_j\}$ for $H^0(X,M^Y)$ that are adapted to the point $p\in Y\cap Z$, in the sense that
$$\mbox{ord}_p(\sigma_{i})=a^{L_Y}_i(p) \mbox{ and } \mbox{ ord}_p(\tau_{j})=a^{M_Y}_j(p),$$
for all relevant $i$ and $j$ (cf. e.g. \cite{EH1}, Lemma 2.3, for the fact that this can be done). Then there are integers $\alpha_i$ and $\beta_j$ defined by
$$\sigma_i^Z=t^{\alpha_i} \sigma_i \mbox{ and } \tau_j^Z=t^{\beta_j}\tau_j.$$
To obtain a formula for the integer $\alpha$ we write $\rho^Y=\sum_{i,j} f_{ij}\sigma_i \otimes \tau_j$, where $f_{ij}\in R$. We have the identity
$$\rho^Z=\sum_{i,j} (t^{\alpha-\alpha_i-\beta_j} f_{ij})(t^{\alpha_i} \sigma_i)\otimes (t^{\beta_j} \tau_j),$$
from which we easily deduce that
$$\alpha=\mbox{max}_{i,j} \{\alpha_i+\beta_j-\nu(f_{ij})\},$$ where $\nu$ denotes the valuation on $R$ (see also \cite{EH4}, Lemma 3.2).

\begin{lemma}\label{sym}
With the above notations, if $\rho_Y:=\rho^Y_{|_Y}\in H^0(Y,L_Y)\otimes H^0(Y,M_Y)$ and $\rho_Z:=\rho^Z_{|_Z}\in H^0(Z,L_Z)\otimes H^0(Z,M_Z)$, then
\begin{center}
$\rm{ord}$$_p(\rho_Y)+\rm{ord}$$_p(\rho_Z)\geq \rm{deg}$$(L_Y)+\rm{deg}$$(M_Y).$
\end{center}
\end{lemma}
\begin{proof} By definition, there exists a pair of indices $(i_1,j_1)$ such that $\nu(f_{i_1j_1})=0$ and
$$\mbox{ord}_p(\rho_Y)=\mbox{ord}_p(\sigma_{i_1,Y})+\mbox{ord}_p(\sigma_{j_1,Y})$$
and clearly $\alpha\geq \alpha_{i_1}+\beta_{j_1}$. To get an estimate on $\mbox{ord}_p(\rho_Z)$ we only have to take into account the pairs of indices $(i,j)$ for which $\alpha_i+\beta_j=\alpha+\nu(f_{ij})\geq \alpha_{i_1}+\beta_{j_1}$. For at least one such pair $(i,j)$ we have that
$$\mbox{ord}_p(\rho_Z)=\mbox{ord}_p(t^{\alpha_i}\sigma_{i, Z})+\mbox{ord}_p(t^{\beta_j}\tau_{j, Z})\geq \alpha_i+\beta_j.$$

\noindent
On the other hand, by applying Lemma \ref{limit} we can write
$$\mbox{ord}_p(\rho_Y)=\mbox{ord}_p(\sigma_{i_{1},Y})+\mbox{ord}_p(\tau_{j_1,Y})\geq \mbox{deg}(L_Y)+\mbox{deg}(M_Y)-\alpha_{i_1}-\beta_{j_1},$$
whence we finally have that
$\mbox{ord}_p(\rho_Z)+\mbox{ord}_p(\rho_Y)\geq \mbox{deg}(L_Y)+\mbox{deg}(M_Y).$
\end{proof}

\section{The class of $\kk$}

We study in detail the Deligne-Mumford compactification of the $K3$ locus.
Since Theorem \ref{equiv} gives four different characterizations for the points of $\KK$, to compute the class of $\kk$ in $\mbox{Pic}(\mm_{10})$ one can try to understand what happens to each of the conditions $(1)-(4)$ as a smooth genus $10$ curve $C$ degenerates to a singular stable curve.

It seems very difficult to understand the degenerations of $[C]\in \KK$ using directly  characterization $(1)$, although it would be highly interesting to have a description of all stable limits of $K3$ sections of genus $g$ (perhaps as curves sitting on Kulikov degenerations of $K3$ surfaces -- cf. also Remark \ref{k3_degenerations}). It also seems almost certain that one cannot use the Wahl map and description $(2)$ to carry out any intersection theoretic computations on $\mm_g$: although the Wahl map can be naturally extended to all stable nodal curves $[C]\in \mm_g$ as $\psi_{\omega_C}:\wedge^2H^0(\omega_C)\rightarrow H^0(\omega_C^{\otimes^2}\otimes K_C)$,
where $\omega_C$ is the dualizing (locally free) sheaf of $C$ and $K_C=\Omega_C^1$, it is easy to see that as soon as $C$ has a disconnecting node
(in particular whenever $C$ is of compact type), $\psi_{\omega_C}$ cannot be surjective for trivial reasons.

Instead, in order to understand $\kk$ we shall use description (4) from Theorem \ref{equiv} together with the set-up developed in \S4.

We  recall some basic things about Brill-Noether divisors on the universal curve $\mm_{g,1}$.
Let us fix positive integers $g,r,d$ and a ramification sequence $\alpha=(0\leq \alpha_0\leq\ldots\leq\alpha_r\leq d-r)$ such that
$$\rho(g,r,d)-w(\alpha)=g-(r+1)(g-d+r)-\sum_{i=0}^r \alpha_i=-1.$$
We define $\cM^r_{g,d}(\alpha)$ to be the locus of pointed curves $(C,p)\in \cM_{g,1}$ such that $C$ carries a $\mathfrak g^r_d$, say $l$, with ramification at $p$ at least $\alpha$, that is, $a_i^l(p)\geq \alpha_i+i$, for $i=0,\ldots,r$. Eisenbud and Harris proved (cf. \cite{EH2}, Theorem 1.2) that the compactification $\mm^r_{g,d}(\alpha)$ is a divisor on $\mm_{g,1}$ which we shall call a Brill-Noether divisor on $\mm_{g,1}$.  They also showed that its class is a linear combination
$$\mm^r_{g,d}(\alpha)\equiv a\cdot BN_1+b\cdot \overline{\mathcal{W}},\mbox{ with }a,b\in \mathbb Q_{\geq 0},$$
where
$$BN_1=(g+3)\lambda-\frac{g+1}{6}\delta_0-\sum_{i=1}^{g-1} i(g-i)\delta_i$$ is the pull-back from $\mm_g$ of the Brill-Noether class and
$\overline{\mathcal{W}}$
is the closure of the Weierstrass locus, considered in \S2 (cf. \cite{EH2}, Theorem 4.1).

By Theorem \ref{equiv} $(4)$ we view $\KK$ as the locus of genus $10$ curves $C$ having a line bundle $M\in W^4_{12}(C)$ such that the multiplication map
\begin{equation}\label{mult}
\mu_M:\mbox{Sym}^2 H^0(M)\rightarrow H^0(M^{\otimes^2}) \mbox{ is not surjective  }(\Leftrightarrow \mbox{ injective}).
\end{equation}
 The class of $\kk$ can be then written then as
$$\kk\equiv A\lambda-B_0\delta_0-B_1\delta_1-\cdots -B_5\delta_5,$$
where $A,B_i\geq 0$. Cukierman and Ulmer showed that $A=7$ and $B_0=1$ (cf. \cite{cukierman-ulmer}, Proposition 3.5), while
Theorem \ref{ineq1} gives the inequalities $B_1\geq 5$ and $B_i\geq 11-i$, for $i=2,\ldots,5$. To compute the coefficients $B_i$ we have to interpret
condition (\ref{mult}) when $C$ is a stable curve of compact type. We  consider the maps $j_i:\mm_{i,1}\rightarrow \mm_{10}$ as in Proposition \ref{attach} and compute the pullback of $\kk$. We have the following result:

\begin{theorem}\label{pullback}
Let $j_i:\mm_{i,1}\rightarrow \mm_{10}$ be the map obtained by attaching a fixed general curve of genus $10-i$. Then for all $1\leq i\leq 4$ the pullback $j^*_i(\kk)$ is a union of Brill-Noether divisors on $\mm_{i,1}$, hence its class is a linear combination of $\overline{\mathcal{W}}$ and  the Brill-Noether class pulled back from $\mm_i$.
\end{theorem}

\begin{remark}
 To make Theorem \ref{pullback} more precise, we can show that $j_1^*(\kk)=0$,  that $j_2^*(\kk)$ is supported on the Weierstrass divisor of $\mm_{2,1}$, while $j_3^*(\kk)$ is supported on the union of the Weierstrass divisor on $\mm_{3,1}$ and the hyperelliptic locus.
\end{remark}

As a corollary, combining this with Proposition \ref{attach}, we get all the coefficients but $B_5$ in the expression of the class of $\kk$, thus proving Theorem \ref{class}.

\noindent
\emph{Proof of Theorem \ref{pullback}.} We shall only describe the case of the map $j_4:\mm_{4,1}\rightarrow \mm_{10}$, the remaining cases involving the maps $j_1,j_2$ and $j_3$ being similar and simpler. We assume that the conclusion of the Theorem already holds for $j_1, j_2$ and $j_3$.
Throughout the proof we consider a genus $10$ curve $X_0=C\cup_p Y$, where $(C,p)$ is a general pointed curve of genus $6$ fixed once and for all, while $(Y,p)$ is an arbitrary smooth pointed curve of genus $4$. We shall prove that we can choose $[C,p]\in \mathcal{M}_{6,1}$ sufficiently general such that for all $[Y,p]\in \mathcal{M}_{4,1}$ outside the union of all Brill-Noether divisors, we have that $[X_0]\notin \kk$. Suppose by contradiction that $[X_0]\in \kk$ and let $\pi:X\rightarrow B$ be a $1$-dimensional family with smooth genus $10$ general fibre $X_t$ sitting on a $K3$ surface and special fibre $X_0'$ semistably equivalent to $X_0$ and obtained from $X_0$ by inserting a  (possibly empty) chain of $\PP^1$'s at the node $p$. According to Theorem \ref{equiv} $(4)$, on a smooth curve $X_t$ near $X_0'$,  there exists a line bundle $L_t\in W^4_{12}(X_t)$ such that the multiplication map
$$\mu_{t}:\mbox{Sym}^2 H^0(X_{t},L_{t})\rightarrow H^0(X_{t},L_{t}^{\otimes^2})$$
is not injective. Take $0\neq \rho_{t}\in \mbox{Ker}(\mu_{t})$ and denote by
$$\{l_Y=\bigl(L_Y,V_Y\subset H^0(L_Y)\bigr), ~l_C=\bigl(L_C, V_C\subset H^0(L_C)\bigr)\}$$
the induced limit $\mathfrak g^4_{12}$ on $X_0$ obtained by restriction from the corresponding limit $\mathfrak g^4_{12}$ on $X_0'$. From general facts about limit linear series we know that there is a $1:1$ correspondence between limit $\mathfrak g^4_{12}$'s on $X_0'$ and $X_0$, and so we may as well assume that $X_0=X_0'$.  According to Section 3 we obtain elements
$$\rho_Y\in \mbox{Ker}\{\mbox{Sym}^2 (V_Y)\rightarrow H^0(L_Y^{\otimes ^2})\}\mbox{ and }\rho_C\in \mbox{Ker}\{\mbox{Sym}^2 (V_C)\rightarrow H^0(L_C^{\otimes ^2})\},$$
such that $\mbox{ord}_p(\rho_Y)+\mbox{ord}_p(\rho_C)\geq 12+12=24$ (cf. Lemma \ref{sym}).

Our assumption that both $(Y,p)$ and $(C,p)$ are Brill-Noether general gives, using the additivity of the Brill-Noether number, that
\begin{equation}\label{bn1}
\rho(6,4,12)-\sum_{i=0}^4(a_i^{l_C}(p)-i)=0,\mbox{ }\mbox{ and}
\end{equation}
\begin{equation}\label{bn2}
\rho(4,4,12)-\sum_{i=0}^4(a_i^{l_Y}(p)-i)=0.
\end{equation}
Since $a_i^{l_C}(p)+a_{4-i}^{l_Y}(p)\geq 12$ for $i=0,\ldots,4$, equalities (\ref{bn1}) and (\ref{bn2}) yield
\begin{equation}\label{bn3}
a_i^{l_C}(p)+a_{4-i}^{l_Y}(p)=12, \mbox{ for }i=0,\ldots,4.
\end{equation}
Moreover, we have the general fact (cf. \cite{EH3}, Proposition 1.2), that a Brill-Noether general pointed curve $(Z,p)\in \cM_{g,1}$ carries a $\mathfrak g^r_d$ with ramification $\geq(\alpha_0,\ldots,\alpha_r)$ at the point $p$ if and only if
\begin{equation}\label{bn4}
g-\sum_{i=0}^r \mbox{max}\{0, \alpha_i+g-d+r\}\geq 0.
\end{equation}
(This is a strengthening of the inequality
\begin{center}
$\rho(g,r,d)-w(p)=\rho(g,r,d)-\sum_{i=0}^r\alpha_i=g-\sum_{i=0}^r(\alpha_i+g-d+r)\geq 0$.)
\end{center}

Conditions (\ref{bn1}), (\ref{bn2}), (\ref{bn3}) and (\ref{bn4}) cut down the number of numerical possibilities for the ramification at $p$ of the limit $\mathfrak g^4_{12}$ on $X_0$ to three. To simplify notations we set $a_i:=a_i^{l_Y}(p)$ and $b_i:=a_i^{l_C}(p)$ for $i=0,\ldots,4$. We have three distinct numerical situations which we shall investigate separately:
\begin{enumerate}
\item $(a_0,a_1,a_2,a_3,a_4)=(4,5,6,9,10)$ and $(b_0,b_1,b_2,b_3,b_4)=(2,3,6,7,8)$
\item $(a_0,a_1,a_2,a_3,a_4)=(4,5,7,8,10)$ and $(b_0,b_1,b_2,b_3,b_4)=(2,4,5,7,8)$
\item $(a_0,a_1,a_2,a_3, a_4)=(4,6,7,8,9)$ and $(b_0,b_1,b_2,b_3,b_4)=(3,4,5,6,8)$.
\end{enumerate}

\noindent
We first study case $(1)$, which will serve as a model for $(2)$ and $(3)$.

Since $\rho_C\in \mbox{Ker}\{\mbox{Sym}^2H^0(C,L_C)\rightarrow H^0(C,L_C^{\otimes ^2})\}$, there must be distinct pairs of indices $(i_1,j_1)\neq (i_2,j_2)$ such that (cf. \S3)
$$\mbox{ord}_p(\rho_C)=b_{i_1}+b_{j_1}=b_{i_2}+b_{j_2}.$$
Similarly, for $Y$ we obtain pairs $(i_1',j_1')\neq (i_2',j_2')$ such that
$$\mbox{ord}_p(\rho_Y)=a_{i_1'}+a_{j_1'}=a_{i_2'}+a_{j_2'}\geq 24-\mbox{ord}_p(\rho_C).$$
Clearly $\mbox{ord}_p(\rho_C)\geq 9$ and case $(1)$ breaks into two subcases:

\noindent \textbf{1$_a$)} $\mbox{ord}_p(\rho_C)=9$, hence $\mbox{ord}_p(\rho_Y)\geq 15(=5+10=6+9).$ This means that $\rho_Y\in \mbox{Sym}^2 H^0\bigl(Y,L_Y(-5p)\bigr)$, where $E:=L_Y(-5p)$ is a $\mathfrak g^3_7$ with vanishing $
(0,1,4,5)$ at $p$. We reach a contradiction by showing that the map $\mbox{Sym}^2H^0(Y,E)\rightarrow H^0(Y,E^{\otimes ^2})$ is injective when $(Y,p)$ is Brill-Noether general. This follows from the base-point-free pencil trick which, applied here, says that the map
$$H^0\bigl(Y,E(-4p)\bigr)\otimes H^0(Y,E)\rightarrow H^0\bigl(Y,E^{\otimes^2}(-4p)\bigr)$$
has kernel $H^0(Y,\OO_Y(4p))$ which is one-dimensional ($p\in Y$ is not a Weierstrass point).

\noindent \textbf{1$_b$)} $\mbox{ord}_p(\rho_C)=10(=3+7=2+8),$ hence $\mbox{ord}_p(\rho_Y)\geq 14(=4+10=5+9).$ This case is more complicated since we cannot reach a contradiction by working with $Y$ alone as we did in $(1_a)$: by counting dimensions it turns out that the map $\mbox{Sym}^2H^0\bigl(Y,L_Y(-4p)\bigr)\rightarrow H^0\bigl(Y, L_Y^{\otimes^2}(-8p)\bigr)$ cannot be injective, so potentially we could  find an element
$\rho_Y\in \mbox{Ker}\{\mbox{Sym}^2(V_Y)\rightarrow H^0(Y,L_Y^{\otimes^2})\}$ with $\mbox{ord}_p(\rho_Y)\geq 14$.

We turn to the genus $6$ curve $C$ instead, and we denote $M:=L_C(-2p)$. Thus $M\in W^4_{10}(C)$ and $a^M(p)=(0,1,4,5,6)$, and we know that there exists an element
$$\gamma \in \mbox{Ker}\{\mbox{Sym}^2H^0(C,M)\rightarrow H^0(C,M^{\otimes ^2})\}$$
such that $\mbox{ord}_p(\gamma)=6(=0+6=1+5)$ (here $\gamma$ is obtained from $\rho_C$ by subtracting the base locus of the linear series $|L_C|$). By degeneration methods we show that such a $\gamma$ cannot exist when $(C,p)\in \cM_{6,1}$ is suitably general.

From general Brill-Noether theory we know that on $C$ there are finitely many line bundles $M\in W^4_{10}(C)$ satisfying $a^M(p)=(0,1,4,5,6)$. They are all of the form $M=K_C\otimes A^{\vee}\otimes \OO_C(4p)$, where $A\in W^1_4(C)$. Since the Hurwitz scheme of coverings $C\stackrel{4:1}\rightarrow \PP^1$ with a genus
$6$ source curve is irreducible, it follows that the variety
$$\{(C,p,M):(C,p)\in \cM_{6,1},M\in W^4_{10}(C)\mbox{ and }a^M(p)\geq (0,1,4,5,6)\}$$
is irreducible as well. Therefore to show that $\gamma$ as above cannot exist when $(C,p)\in \cM_{6,1}$ is general, it will be enough to prove the following:

\noindent
\textbf{Claim:} The general genus $6$ pointed curve $(C,p)$ with $C\subset \PP^4, \mbox{deg}(C)=10$ and $a(p)=a^{\OO_C(1)}(p)=(0,1,4,5,6)$ does not sit on a quadric $Q\subset \PP^4$ with $\mbox{ord}_p(Q)\geq 6$.

\noindent
\emph{Proof of claim.}  By semicontinuity, it is enough to construct one single such embedded curve. We start with the smooth monomial curve $\Gamma=\nu(\PP^1)$, where $\nu:\PP^1\rightarrow \PP^4, \nu(t)=[1:t:t^4:t^5:t^6]$. We set $p:=\nu(0)\in \Gamma$ and then $a(p)=(0,1,4,5,6)$. If $x_0,x_1,x_2,x_3$ and $x_4$ are the coordinates on $\PP^4$ then $\Gamma$ is contained in three quadrics:
$$Q_1:x_1x_2-x_0x_3=0, \mbox{ }Q_2:x_0x_4-x_1x_3=0\mbox{ and } Q_3:x_2x_4-x_3^2=0.$$
Note that $\mbox{ord}_p(Q_1)=5$ and $\mbox{ord}_p(Q_2),\mbox{ord}_p(Q_3)\geq 6$.
Take now a general hyperplane $H\subset \PP^4$ and denote by $\{p_1,\ldots,p_6\}=H\cap \Gamma.$ The points $p_1,\ldots, p_6$ will be in general position in $H=\PP^3$ and let $R\subset \PP^3$ be the unique twisted cubic passing through them. We choose an elliptic quartic $E\subset \PP^3$ that
passes through $p_1,\ldots,p_6$ and set $C:=\Gamma \cup E$. Thus $C\subset \PP^4$ is a stable curve of genus $6$ with $\mbox{deg}(C)=10$ and vanishing $a(p)=(0,1,4,5,6)$ at the smooth point $p$.  To show that $C$ can be smoothed in $\PP^4$ while preserving the ramification at $p$ it is enough to notice that $A:=\omega_C(4p)\otimes \OO_C(-1)$ is a $\mathfrak g^1_4$ on $C$ and use that every pencil on a stable curve can be smoothed to nearby smooth curves (cf. e.g. \cite{EH1}). Alternatively one can prove this smoothing result by employing the methods from \cite{hh}.

The quadrics in $\PP^4$ containing $\Gamma$ with order $\geq 6$ at $p$ will be those in the pencil $\{\lambda Q_2+\nu Q_3\}_{[\lambda:\nu]\in \PP^1}$. We prove that no quadric in this pencil can contain a general elliptic quartic $E\subset H$ that passes through $p_1,\ldots, p_6$, which will settle the claim. This follows from a dimension count: the space of elliptic quartics $E\subset H$ through $p_1,\ldots,p_6$ is $4$-dimensional, while $\{B_{[\lambda:\nu]}=H\cap (\lambda Q_2+\mu Q_3=0)\}_{[\lambda:\nu] \in \PP^1}$ is a pencil of quadrics in $H$ containing the twisted cubic $R$. It is immediate to see that there are only $\infty^2$ elliptic quartics on each of the quadrics $B_{[\lambda:\nu]}$ in the pencil, which finishes the proof of the claim and finally takes care of case $(1_b)$.

\noindent \textbf{1$_c$)} $\mbox{ord}_p(\rho_C)\geq 11$, in which case $\mbox{ord}_p(\rho_C)=14(=7+7=6+8)$, and this is impossible since the map $\mbox{Sym}^2 H^0\bigl(C,L_C(-6p)\bigr)\rightarrow H^0\bigl(C,L_C^{\otimes^2}(-12p)\bigr)$ is injective.

Cases $(2)$ and $(3)$ are similar to $(1)$. For instance $(3)$ breaks into the following subcases:

\noindent \textbf{3$_a$)} $\mbox{ord}_p(\rho_C)=8(=3+5=4+4)$, which can be dismissed right away by looking at the ramification on $Y$ (it would imply that $\mbox{ord}_p(\rho_Y)\geq 16(=7+9=8+8)$ which is impossible).

\noindent \textbf{3$_b$)} $\mbox{ord}_p(\rho_C)=9(=3+6=4+5)$, when also $\mbox{ord}_p(\rho_Y)=15(=9+6=7+8)$. This case is ruled out by applying the base-point-free pencil trick on $Y$ just like we did in case $(1_a)$.

\noindent \textbf{3$_c$)} $\mbox{ord}_p(\rho_C)\geq 10$. This case is similar to $(1_b)$. We set $M:=L_C(-3p)\in W^4_9(C)$.
Then we have the vanishing sequence $a^M(p)=(0,1,2,3,5)$ and there exists an element
$\gamma \in \mbox{Ker}\{\mbox{Sym}^2 H^0(C,M)\rightarrow H^0(C,M^{\otimes ^2})\}
$
 with $\mbox{ord}_p(\gamma)\geq 4$.
\newline
\indent
Riemann-Roch gives that $M=K_C\otimes A^{\vee}\otimes \OO_C(5p)$, where $A\in W^1_6(C)$ is such that $|A|(-5p)\neq \emptyset$. Since the Hurwitz scheme of coverings $C\stackrel{6:1}\rightarrow \PP^1$ with a genus $6$ source curve and a ramification point $p$ of index $5$ is irreducible, again, it will be sufficient to construct a single smooth curve $C\subset \PP^4$ of genus $6$, $\mbox{deg}(C)=9$ and having a point $p\in C$ with $a(p)=(0,1,2,3,5)$ and such that $C$ does not lie on a quadric $Q\subset \PP^4$ with $\mbox{ord}_p(Q)\geq 4$.

To construct $C$, we start with a pointed elliptic curve $(\Gamma,p)$ which we embed in $\PP^4$ by the linear series $|5p|$. We choose a general hyperplane $H\subset \PP^4$ and denote by $\{p_1,\ldots,p_5\}=\Gamma \cap H$. Then we fix a general elliptic quartic $E\subset H$ through $p_1,\ldots,p_5$ and set $C:=\Gamma \cup E$. We have that $C\subset \PP^4$ is a stable curve of genus $6$ and degree $9$ with vanishing $a(p)=(0,1,2,3,5)$. There is a $5$-dimensional family of quadrics in $\PP^4$ containing $\Gamma$ and a $3$-dimensional subfamily of quadrics $Q\subset \PP^4$ with $\mbox{ord}_p(\Gamma)\geq 4$ (cf. \cite{Hu}, Proposition IV.2.1). A dimension count similar to the one in $(1_b)$ establishes that none of these $\infty^3$ quadrics will contain a general elliptic quartic $E\subset H$ that passes through $p_1,\ldots,p_5$.

Thus we have proved that there is $[C,p]\in \mathcal{M}_{6,1}$ such that for any $[Y,p]\in \mathcal{M}_{4,1}$ which is outside the Brill-Noether divisors, we have that $[Y\cup_p C]\notin \kk$. The same conclusion holds if $[Y,p]$ is a
general element of the divisor $\Delta_0$ on $\mm_{4,1}$ (the base-point-free pencil trick and the usual Brill-Noether theory used in the proof carries through to general irreducible nodal curves, (cf. \cite{EH3}, Theorem 1.1)).
Therefore $j_4^*(\kk)$ is a union of Brill-Noether divisors and (possibly)
multiples of the boundary divisors $\Delta_1, \Delta_2$ and $\Delta_3$, hence we can write $j_4^*(\kk)\equiv \alpha \overline{\mathcal{W}}+\beta BN_1+c_1 \Delta_1+
c_2 \Delta_2+c_3\Delta_3$, for $\alpha, \beta, c_1, c_2, c_3 \in \mathbb Q_{\geq 0}$.

Since we assume having already proved that $j_i^*(\kk)$ is a union of Brill-Noether divisors for $i=1,2,3$, we have that
$$\kk\equiv 7\lambda-\delta-5\delta_1-9\delta_2-12\delta_3-B_4\delta_4-B_5\delta_5.$$
By identification we obtain $\alpha=7/5, \beta=6/5, B_4=14$ and $c_1=c_2=c_3=0$, which finishes the proof. Note that we have also proved Theorem \ref{class} in the process.
\hfill $\Box$

\begin{remark} An easy calculation shows that for $6\leq i\leq 9$, $j_i^*(\kk)$ is \emph{not} a combination of Brill-Noether divisors on $\mm_{i,1}$, a fact that can be also seen in the last part of the proof of Theorem \ref{pullback}.
We have been unable to determine $j_5^*(\kk)$. \end{remark}

\section{The Slope Conjecture and the Iitaka dimension of the Brill-Noether linear
system on $\mm_g$}

This section, which should be considered at least partially joint with Sean Keel, emphasizes another
somewhat surprising use of curves on $K3$ surfaces, this time in genus
$11$. We note, based
on results of Mukai, that on $\mm_{11}$ there exist \lq \lq many strongly independent"
effective divisors of minimal slope. This seems to contradict  earlier beliefs; see below for a more
precise formulation.
We denote by $\overline{\mathcal{F}}_g$ the Baily-Borel
compactification of $\mathcal{F}_g$ (see e.g. \cite{Loo} for a general reference).
For a $\mathbb Q$-Cartier divisor $D$ on a variety $X$ we denote by
$\kappa(D)=\kappa(X,D)$ its Iitaka dimension.

Recall that Harris and Morrison have conjectured that $s(D)\geq 6+12/(g+1)$ for any effective divisor $D$ on $\mm_g$, with equality when $g+1$ is composite.
In this case, the quantity $6+12/(g+1)$ is the
slope of the Brill-Noether divisors on $\mm_g$. Harris and Morrison
wondered whether  all
effective divisors of slope $6+12/(g+1)$ should consist of curves having some special character
(cf. \cite{harris-morrison}, p. 324).
In this direction, they proved that on $\mm_3$ the only irreducible
divisor of slope $9\bigl(=6+12/(g+1)\bigr)$ is the hyperelliptic locus
$\mm_{3,2}^1$, and that on $\mm_5$ the only irreducible divisor of
slope $8$ is the trigonal locus $\mm_{5,3}^1$. Moreover, a standard
argument involving pencils of plane curves shows that on $\mm_8$, the
Brill-Noether divisor $\mm_{8,7}^2$ is the only irreducible divisor of
slope $22/3$, while a slightly more involved argument using pencils of $4$-gonal curves proves a similar conclusion for the $4$-gonal divisor $\mm_{7,4}^1$ on $\mm_7$. One way to rephrase these results is to say that for
$g=3,5, 7, 8$,  every effective
$\QQ$-divisor $D$ with $s(D)=6+12/(g+1)$ has Iitaka dimension $\kappa (D)=0$, and each such $D$ is an effective combination of the (unique) Brill-Noether divisor $\mm_{g,d}^r$ and boundary divisors $\Delta_1,\ldots, \Delta_{[g/2]}.$ In the case when $g+1$ is multiply divisible and there are several Brill-Noether divisors $\mm_{g,d}^r$, it is natural to ask the following:

\begin{question}\label{strongslope}
For composite $g+1$, if $a$ denotes the number of distinct Brill-Noether divisors on $\mm_g$, is it true that $\kappa(\mm_g,D)\leq a-1$ for any effective $\mathbb Q$-divisor $D$ on $\mm_g$ with slope $s(D)=6+12/(g+1)$?
\end{question}

The first interesting case is $g=11$, when
there are two distinct Brill-Noether divisors $\mm_{11,6}^1$ and $\mm_{11,9}^2$.
We are asking  whether $\kappa(BN+\sum_{i=1}^5 a_i\delta_i)=1$ for all $a_i\geq 0, i=1,\ldots,5$, where $BN:=7\lambda-\delta_0-5\delta_1-9\delta_2-12\delta_3-14\delta_4-15\delta_5$ is the Brill-Noether class.  (The next interesting case is that of $\mm_{23}$. This case was studied
extensively in \cite{F}, and the Iitaka dimension of the Brill-Noether system is conjectured to
give a positive answer to the question above.)
We show that at least  in this genus, the answer to Question \ref{strongslope} is far from being positive:

\begin{proposition}\label{big_iitaka}
There exist effective divisors on $\mm_{11}$ of slope $7$ and having
Iitaka dimension equal to $19$.
For example $\kappa(\mm_{11}, BN+4\delta_3+7\delta_4+8\delta_5)= 19$.
\end{proposition}
\begin{proof} By work of Mukai (cf. \cite{Mukai2}) there exists a rational map
$\phi:\mm_{11} - ->\overline{\mathcal{F}}_{11}$ which sends a general curve $[C]\in \cM_{11}$ to $[S,C]\in \mathcal{F}_{11}$, where $S$ is the unique $K3$ surface containing $C$.
We denote by $A$ the indeterminacy locus of $\phi$. Note that since $\mm_{11}$ is normal, $\mbox{codim}(A,\mm_{11})\geq 2$. For a $\mathbb Q$-Cartier divisor $D$ on $\ff_{11}$ we define the pull-back $\phi^*D$, for example as $p_*q^*(D)$, where $p:\Sigma\rightarrow \mm_{11}$ and $q:\Sigma\rightarrow \ff_{11}$ are the projections from the closure of the graph of $\phi$.
It is easy to check that $\kappa(\mm_{11},\phi^*D)\geq \kappa(\ff_{11},D)$.

Take $D$ now to be any ample effective divisor on $\ff_{11}$ (e.g.  the zero locus of an automorphic form on the period space of $K3$ surfaces). Let us write $\phi^*(D)\equiv a\lambda-\sum_{i=0}^5 b_i\delta_i$, and we claim that $a/b_0=7$. We choose a general $[S,C]\in \mathcal{F}_{11}$ such that $\mbox{Pic}(S)=\mathbb Z C$ and $S$ is the only $K3$ surface containing $C$.
We also pick a Lefschetz pencil on $S$, giving rise to a curve $B$ in $\mm_{11}$. Since $B$ fills-up $\mm_{11}$, we can assume that $B\cap A=\emptyset$
and that $\phi(b)=[S,C]$ for a general $b\in B$, hence $\phi(B)=\{[S,C]\}$, that is, $\phi$ contracts the curve $B$. Then $B\cdot \phi^*(D)=0$. On the other hand we have that $B\cdot \lambda=g+1=12, B\cdot \delta_0=6(g+3)=84$ and $B\cdot \delta_i=0$ for $i\geq 1$ (cf. Lemma \ref{nos1}), thus we obtain that $a/b_0=B\cdot \delta_0/B\cdot \lambda =7$.

We decompose $\phi^*(D)$ into its moving and fixed part, $\phi^*(D)=M+F$, and choose a general member $D'\in |M|$. Then $D'\equiv a'\lambda-\sum_{i=0}^5 b_i'\delta_i$ contains none of the boundaries $\Delta_1,\ldots, \Delta_5$ as components, and from Theorem \ref{ineq1} we obtain that $a'/b_0'=a/b_0=7, b_1'\geq 5b_0'$ and $b_i'\geq (11-i)b_0'$ for $i=2,\ldots ,5$. If we set $E:=BN+4\delta_3+7\delta_4+8\delta_5$, it follows then that $\kappa(\mm_{11},E)\geq \kappa(\mm_{11},D')= \kappa(\mm_{11},\phi^*(D))\geq 19$.

On the other hand we claim that $\kappa(\mm_{11},E)\leq 19$. Indeed, since $B\cdot E=0$ we have that the rational map associated to any multiple of $E$ contracts the $11$-dimensional family of curves corresponding to the linear system $|C|$ on $S$, thus $\kappa(\mm_{11},E)\leq \mbox{dim}(\mm_{11})-11=19$.
\end{proof}

\section{Further applications}

\subsection*{\bf{The $K3$ locus vs. Brill-Noether loci.}}
For positive integers $g,r,d$ such that $\rho(g,r,d)<0$, a  Brill-Noether locus $\cM_{g,d}^r$ on $\cM_g$ is defined as the locus of curves carrying a $\mathfrak g^r_d$, and we denote by $\mm_{g,d}^r$ its Deligne-Mumford compactification. When $\rho(g,r,d)=-1$, it is known that $\cM_{g,d}^r$ is an irreducible divisor (cf. \cite{EH2}). Virtually all of our knowledge about the effective cone of $\mm_g$ comes from the study of the Brill-Noether divisors $\mm_{g,d}^r$. To emphasize the significant difference between $\KK$ and $\cM_{g,d}^r$ as well as the potential for getting new information on the birational geometry of $\mm_g$ by studying higher genus analogues/generalizations of $\KK$, we shall compare their behavior under the flag map $\phi:\mm_{0,g}\rightarrow \mm_g$ obtained by attaching a fixed elliptic curve to the marked points of a $g$-pointed rational curve. A key result of modern Brill-Noether theory is the following:

\begin{theorem}(\cite{EH3}, Theorem 1.1)
Whenever $\rho(g,r,d)<0$ we have that $\phi^*(\mm_{g,d}^r)=0$, that is, a flag curve of genus $g$ with $g$ elliptic tails is Brill-Noether general.
\end{theorem}
This result, besides offering a proof of the Brill-Noether Theorem,  can be employed to compute almost all the coefficients of the class of $\mm_{g,d}^r$ when $\rho(g,r,d)=-1$ (cf. \cite{EH3}, Theorem 3.1). In contrast, for the locus of  $K3$ sections $\KK_g$ we have the following result:
\begin{proposition}\label{tails}
The flag locus of genus $g$ is entirely contained in $\kk_g$, that is, $\phi(\mm_{0,g})\subset \kk_g$.
\end{proposition}
\begin{proof} Let us denote by $E$ the elliptic tail which appears in the definition of $\phi$ and denote by $S$ the ($K3$) Kummer surface associated with $E\times E$. Pick $a\in E[2]$ an element of order $2$ and denote by $R$ the strict transform of $\{a\}\times E$. Then $R$ is a smooth rational curve and $R^2=-2$. Pick also $x_1,\ldots,x_g\in E-E[2]$ arbitrary points and denote by $C_i$ the strict transform of $E\times \{x_i\}$ for $i=1,\ldots,g$. Then all $C_i$ are elliptic curves isomorphic to $E$ such that $C_i\cdot R=1$, hence
$R+C_1+\cdots+C_g$ is a flag curve of genus $g$ sitting on the $K3$ surface $S$.
\end{proof}

Theorem \ref{class} gives detailed information about \emph{degenerate K3 sections}, that is, stable curves that are limits of smooth $K3$ sections.
For instance one can prove the following:

\begin{proposition}
Every stable curve of genus $10$ with five tails of genus two is a degenerate $K3$ section. Every genus $10$ curve with one elliptic and three genus three tails is a degenerate $K3$ section.
\end{proposition}
\begin{proof}
We only consider the first case. We look at the map $m:\mm_{0,5}\rightarrow \mm_{10}$ obtained by attaching five arbitrary genus two tails at the marked points $x_1,\ldots ,x_5$ of each element from $\mm_{0,5}$. If $B_2$ denotes the boundary divisor of $\mm_{0,5}$ corresponding to singular rational curves of type $(2,3)$, then it is easy to see
that
$$m^*(\delta_4)=B_2,m^*(\delta_2)=-\sum_{i=1}^5 \psi_{x_i}=-\frac{3}{2}B_2, m^*(\lambda)=m^*(\delta_i)=0 \mbox{ for }i\in \{0,1,3,5\}.$$
Thus $m^*(\kk)=-\frac{3}{2}B_2$ and so $\mbox{Im}(m)\subset \kk$.
\end{proof}

\begin{remark}\label{k3_degenerations}
It would be interesting to realize explicitly such genus $10$ curves with five genus two tails as
sections of some $K3$ surfaces. Such a surface will
necessarily be a degenerate one. In fact, if $C$ is a singular genus $g$ curve of compact type sitting on a smooth $K3$ surface $S$, then using well known facts about linear systems on $K3$ surfaces, one can prove that  $S$ is elliptic and $C$ consists only of rational and elliptic curves.
\end{remark}

\subsection*{\bf The Kodaira dimension of $\mm_{10,n }$.}
It is known that for each $g\geq 3$, there is an
integer $f(g)$ such that $\mm_{g,n}$ is of general type for all $n\geq f(g)$
(cf. \cite{Lo}, Theorem 2.4). For those values of $g$ for which
$\mathcal{M}_g$ is unirational (or more generally
$\kappa(\mm_g)=-\infty$), a natural question to ask is to determine
$f(g)$. We show that Theorem \ref{class} can be used to
give an answer to this question for $g=10$.
Recall first that one has the formula for the canonical class (cf. \cite{Lo}, Theorem 2.6)
$$K_{\mm_{g,n}}\equiv 13\lambda-2\delta_{irr}+\sum_{i=1}^n \psi_i-2\sum_{i\geq 0,S}\delta_{i:S}-\sum_S\delta_{1:S}.$$
To prove that $K_{\mm_{g,n}}$ is effective for certain $g$ and $n$ we are going to use besides the divisor $\kk$, the effective divisor $D$ on $\mm_{g,g}$ consisting of pointed genus $g$ curves $(C,p_1,\ldots, p_g)$ such that $h^0(C,p_1+\cdots +p_g)\geq 2$. The class of $D$ has been computed in \cite{Lo}, Theorem 5.4, and we have the formula
$$D\equiv -\lambda+\sum_{i=1}^g \psi_i-\sum_{i\geq 0, S} c_{i:S}\  \delta_{i:S},$$
where $c_{0:S}\geq 2$ and $c_{i:S}> 0$, for all relevant $i$ and $S$ (note that the coefficient of $\delta_{irr}$ is $0$).

\begin{proposition}\label{gentype}
The Kodaira dimension of $\mm_{10,10}$ is $\geq 0$. For all $n\geq 11$ we have that  $\mm_{10,n}$ is of general type.
\end{proposition}
\begin{proof} We fix an integer $n\geq 10$ and denote by $\pi_n:\mm_{10,n}\rightarrow \mm_{10}$ the morphism forgetting all the marked points. We consider two effective divisors on $\mm_{10,n}$: firstly, the pullback of the $K3$ locus
$$\mathfrak{K}_n:=\pi_n^*(\kk)\equiv 7\lambda-\delta_{irr}-5\sum_{S} \delta_{1:S}-9\sum_S \delta_{2:S}-12\sum_S \delta_{3:S}-14\sum_S \delta_{4:S}-B_5\sum_S \delta_{5:S}.$$
Secondly, if for any $S\subset \{1,\ldots,n\}$ with $|S|=10$, we denote by $\pi_S:\mm_{10,n}\rightarrow \mm_{10,10}$ the morphism forgetting the marked points from $S^c$, by averaging the pullbacks of $D$, we define the effective divisor on $\mm_{10,n}$
$$\mathfrak{D}_n:=\frac{1}{{n\choose 10}} \sum_{|S|=10} \pi_S^*(D)\equiv -\lambda+\frac{10}{n}\sum_{i=1}^n \psi_i-\sum_{i\geq 0,S}\tilde{c}_{i:S}\ \delta_{i:S},$$
where $\tilde{c}_{0:S}\geq 2$ and $\tilde{c}_{i:S}>0$. A check yields that  $K_{\mm_{10,10}}-2\mathfrak{K}_{10}-\mathfrak{D}_{10}$ is an effective combination of boundary classes $\delta_{i:S}$, thus proving that $\kappa(\mm_{10,10})\geq 0$. Similarly, for $n\geq 11$, one can find positive constants $\alpha, \beta, a,b_{irr}$ and $c_i$ for $1\leq i\leq n$, with $a/b_{irr}$ arbitrarily large,  such that
$K_{\mm_{10,n}}-\alpha \mathfrak{K}_n-\beta
\mathfrak{D}_n-(a\lambda-b_{irr}\delta_{irr}+\sum_{i=1}^n c_i\psi_i)$
is a positive  combination of boundary classes $\delta_{i:S}$.  As
the class
$a\lambda-b_{irr}\delta_{irr}+\sum_{i=1}^n c_i\psi_i$ is big for such a
choice (cf. e.g. \cite{Lo}, Theorem 2.9), this implies that $\mm_{10,n}$ is of general type.
\end{proof}

\noindent
Using again the divisor $\kk$ we can show that Proposition \ref{gentype} is optimal:

\begin{proposition}
The Kodaira dimension of $\mm_{10,n}$ is $-\infty$ for $n\leq 9$.
\end{proposition}
\begin{proof}
We are only going to prove this for $n=9$, which will
imply the same conclusion for lower $n$. We consider the divisor
$\mathfrak{K}_9=\pi_9^*(\kk)$ on $\mm_{10,9}$ and a general point
$[C,p_1,\ldots, p_9]\in \mathfrak{K}_9$
corresponding to a curve $C$ sitting on a $K3$ surface $S$. Since
$\mbox{dim}|\mathcal{O}_S(C)|=10$, from the generality of $C\subset S$
and of the points $p_i\in C$, it follows that
$|\OO_S(C)\otimes \I_{\{p_1,\ldots, p_9}\}|$ is a pencil giving rise to a curve $R\subset \mm_{10,9}$. One
finds that $R\cdot \lambda=11, R\cdot \delta_{irr}=78, R\cdot
\psi_i=1$ for all
$1\leq i\leq 9$ while $R\cdot \delta_{i:S}=0$ for all $i$ and $S$
(cf. Lemma \ref{nos2} and completely similar calculations). It follows that
$R\cdot K_{\mm_{10,9}}=-4$.

Next we define a second curve $T\subset \Delta_{irr}\subset
\mm_{10,9}$, obtained from a general pointed curve
$[Y,p,p_1,\ldots, p_9]\in \mm_{9,10}$ by identifying the fixed point $p$ with a moving point $y\in Y$.
A standard calculation shows that $T\cdot \lambda=0, T\cdot
\delta_{irr}=-2g(Y)=-18, T\cdot \delta_{1:\emptyset}=1$
and $T\cdot \psi_i=1$ for $1\leq i\leq 9$, while $T$ stays away from all other boundary divisors on $\mm_{10,9}$.

We assume that $K_{\mm_{10,9}}$ is effective and write
$K_{\mm_{10,9}}\equiv m\cdot \mathfrak{K}_9+n\cdot \Delta_{irr}+E$ for
$m,n\in \mathbb Z_{\geq 0}$, where
$E$ is an effective divisor that contains neither $\mathfrak{K}_9$ nor
$\Delta_{irr}$. Since the curves of the type of $R$ fill-up the
divisor $\mathfrak{K}_9$, while those of the type of
$T$ fill-up $\Delta_{irr}$, we must have that $R\cdot E\geq 0$ and
$T\cdot E\geq 0$. A direct check shows that
these inequalities are not compatible with the conditions $m,n\geq 0$.
\end{proof}

\end{document}